\documentclass[11pt]{article}

\usepackage{amsfonts,amsmath,latexsym,color,epsfig,hyperref,enumitem}
\newtheorem{theorem}{Theorem}
\newtheorem{lemma}{Lemma}
\newtheorem{proposition}{Proposition}
\newtheorem{conjecture}{Conjecture}

\def\qed{\ifvmode\mbox{ }\else\unskip\fi\hskip 1em plus 10fill$\Box$}
\newenvironment{proofof}[1]
      {\medskip\noindent{\bf Proof of #1:}\hspace{1mm}}
      {\hfill$\Box$\medskip}
\input{epsf}

\makeatletter
\def\Ddots{\mathinner{\mkern1mu\raise\p@
\vbox{\kern7\p@\hbox{.}}\mkern2mu
\raise4\p@\hbox{.}\mkern2mu\raise7\p@\hbox{.}\mkern1mu}}
\makeatother

\title{\vspace{-0.7cm}On the Zarankiewicz problem for graphs with bounded VC-dimension}
\author{Oliver Janzer\thanks{Department of Mathematics, ETH Z\"urich, Switzerland. Email: {\tt oliver.janzer@math.ethz.ch}. Research supported by an ETH Z\"urich Postdoctoral Fellowship 20-1 FEL-35.} \and Cosmin Pohoata\thanks{Department of Mathematics, Yale University, USA. Email: {\tt andrei.pohoata@yale.edu}.}}
\date{}

\begin{document}
\maketitle

\begin{abstract}
The problem of Zarankiewicz asks for the maximum number of edges in a bipartite graph on $n$ vertices which does not contain the complete bipartite graph $K_{k,k}$ as a subgraph. A classical theorem due to K\H{o}v\'ari, S\'os, and Tur\'an says that this number of edges is $O\left(n^{2 - 1/k}\right)$. An important variant of this problem is the analogous question in bipartite graphs with VC-dimension at most $d$, where $d$ is a fixed integer such that $k \geq d \geq 2$. A remarkable result of Fox, Pach, Sheffer, Suk, and Zahl [{\textit{J. Eur. Math. Soc. (JEMS)}}, no. 19, 1785--1810] with multiple applications in incidence geometry shows that, under this additional hypothesis, the number of edges in a bipartite graph on $n$ vertices and with no copy of $K_{k,k}$ as a subgraph must be $O\left(n^{2 - 1/d}\right)$. This theorem is sharp when $k=d=2$, because by design any $K_{2,2}$-free graph automatically has VC-dimension at most $2$, and there are well-known examples of such graphs with $\Omega\left(n^{3/2}\right)$ edges. However, it turns out this phenomenon no longer carries through for any larger $d$. 

We show the following improved result: the maximum number of edges in bipartite graphs with no copies of $K_{k,k}$ and VC-dimension at most $d$ is $o(n^{2-1/d})$, for every $k \geq d \geq 3$.

\end{abstract}

\section{Introduction}

The problem of Zarankiewicz is a central problem in extremal graph theory. It asks for the maximum number of edges $\operatorname{ex}(n,K_{k,k})$ in a bipartite graph on $n$ vertices, where each side of the bipartition contains $n/2$ vertices and which does not contain the complete bipartite graph $K_{k,k}$ as a subgraph. In 1954, K\H{o}v\'ari, S\'os, and Tur\'an \cite{KST54} proved that this number of edges is at most $c_{k}n^{2 - \frac{1}{k}}$, for some positive constant $c_{k}$ which depends only on $k$. Classical constructions of Reiman and Brown show that this bound is tight for $k=2,3$ (see \cite{PA95}).
However, the Zarankiewicz problem for $k \geq 4$ remains one of the most challenging unsolved
problems in combinatorics.  The best lower bound for $k=4$ simply comes from the Brown construction \cite{Br66}, namely 
$$\operatorname{ex}(n,K_{4,4}) \geq \operatorname{ex}(n,K_{3,3}) = \Omega \left(n^{5/3}\right).$$
For $k=5,6$, the best lower bounds are due to Ball and Pepe \cite{BP12} and come from the norm graph construction of Alon, R\'onyai and Szab\'o \cite{ARS99} (originally considered for the asymmetric Zarankiewicz problem regarding bipartite graphs containing no copies of $K_{4,7}$), i.e.
$$\operatorname{ex}(n,K_{6,6}) \geq \operatorname{ex}(n,K_{5,5}) = \Omega(n^{7/4}).$$
For $k \geq 7$, the best construction comes from a result of Bohman and Keevash \cite{BK10} on random graph processes, and is of the form
$$\operatorname{ex}(n,K_{k,k}) = \Omega\left(n^{2 - 2/(k+1)} (\log n)^{1/(k^2-1)}\right).$$

An important variant of this problem is the analogous question in bipartite graphs with {\it{VC-dimension}} at most $d$, where $d$ is a fixed integer not larger than $k$. The {\it{VC-dimension of a set system}} $\mathcal{F}$ on the ground set $V$ is the largest integer $d$ for which there exists a $d$-element set $S \subset V$ such that for every subset $B \subset S$, one can find a member $A \in \mathcal{F}$ with $A \cap S = B$. The VC-dimension, introduced by Vapnik and Chervonenkis \cite{VC71}, is one of the most useful combinatorial parameters that measures the complexity of graphs and hypergraphs. Over the years it has proven tremendously important in many areas in (and out of) combinatorics and mathematics, in general. We refer to \cite{FPS20} for a nice discussion (and a further remarkable application). In order to define the VC-dimension of a bipartite graph $G = (A,B,E)$ with vertex set $A \cup B$ and edge set $E \subset A \times B$, the standard convention is to make a choice between $A$ and $B$, and then define the {\it{VC-dimension of $G = (A,B,E)$ with respect to $A$}} to be the VC-dimension of the set system of neighborhoods of vertices $b \in B$ (regarded as subsets of $A$), and vice versa (see also \cite{FPSSZ17}). In this paper, we will always choose the first side of the bipartition $A$ as the ground set, and so we shall say that the {\it{VC-dimension of a bipartite $G = (A,B,E)$}} is the VC-dimension of $G$ with respect to $A$. 

With this terminology, we now state the following remarkable result of Fox, Pach, Sheffer, Suk and Zahl from \cite[Theorem 2.1]{FPSSZ17}, which goes below the K\H{o}v\'ari, S\'os, and Tur\'an upper bound for the Zarankiewicz problem, in the presence of bounded VC-dimension.

\begin{theorem}\label{FPSSZ1}
Let $k$ and $d$ be integers such that $k \geq d \geq 2$. Let $G = (A,B,E)$ be a bipartite graph with $|A|=|B|=n/2$ with VC-dimension at most $d$. Then, if $G$ is $K_{k,k}$-free, we have
$$|E(G)| \leq c n^{2-1/d},$$
for some positive constant $c=c(d,k)$. 
\end{theorem}

Using a slightly stronger version of Theorem \ref{FPSSZ1} for bipartite graphs $G=(A,B,E)$ with $|A|$ and $|B|$ of possibly different sizes and with {\it{dual shatter function of degree at most d}}, the observation that semi-algebraic bipartite graphs of {\it{bounded description complexity}} also have bounded VC-dimension (a consequence of the Milnor-Thom theorem), and a standard amplification trick via the so-called {\it{cutting lemma}} of Chazelle \cite{Ch93}, the authors deduced the following general result.

\begin{theorem} \label{FPSSZ2}
Let $G = (P,Q,E)$ be a semi-algebraic bipartite graph in $\left(\mathbb{R}^{d_{1}},\mathbb{R}^{d_{2}}\right)$ such that $E$ has description complexity at most $t$, $|P|=m$, and $|Q|=n$. If $G$ is $K_{k,k}$-free, then
$$|E(G)| \leq c_{1}\left((mn)^{2/3}+m+n\right)$$
for $d_{1}=d_{2}=2$, and 
$$|E(G)| \leq c_{2} \left(m^{\frac{d_{2}(d_{1}-1)}{d_{1}d_{2}-1} + \epsilon} n^{\frac{d_{1}(d_{2}-1)}{d_{1}d_{2}-1}} + m + n \right)$$
for every $\epsilon > 0$, for all $d_{1},d_{2}$. Here $c_{1}=c_{1}(t,k)$ and $c_{2}=c_{2}(d_{1},d_{2},t,k,\epsilon)$.
\end{theorem}

Since in this paper we will not be discussing anything specific to semi-algebraic graphs, we won't attempt to make Theorem \ref{FPSSZ2} (and its connection to Theorem \ref{FPSSZ1}) more precise by providing the required definitions, but we invite the interested reader to consult \cite{FPSSZ17}. Nevertheless, it is worth mentioning at least that this theorem has multiple applications in combinatorial geometry, as the case $d_{1}=d_{2}=2$ already directly implies the celebrated Szemeredi-Trotter theorem \cite{ST83}.

The main result of this paper concerns Theorem \ref{FPSSZ1}. Going back to the statement, it is important to note that it is sharp when $k=d=2$. This is because, by design, every $K_{2,2}$-free bipartite graph has $VC$-dimension at most $2$, and so the well-known examples of such graphs with $\Omega\left(n^{3/2}\right)$ edges automatically serve as constructions which match the upper bound from Theorem \ref{FPSSZ1}. However, it turns out this phenomenon does not extend for larger $k$. For example, already when $k=3$, it is not difficult to check that the dense $K_{3,3}$-free graph from Brown's construction \cite{Br66} has VC-dimension equal to $4$, so it does not anymore provide a matching lower bound. This suggests that a potential improvement of Theorem \ref{FPSSZ1} might be possible when $k \geq d \geq 3$. 

In what follows, we confirm this suspicion by proving the following improved result.

\begin{theorem}\label{main}
Let $k$ and $d$ be integers such that $k \geq d \geq 3$. Let $G = (A,B,E)$ be a bipartite graph with $|A|=|B|=n/2$ with VC-dimension at most $d$. Then, if $G$ is $K_{k,k}$-free, we have
$$|E(G)| = o_{n \to \infty}\left(n^{2-1/d}\right).$$
\end{theorem}

The proof of Theorem 3, which we will discuss in Section 2, is very different than the proof of Theorem 1 from \cite{FPSSZ17}, since our method does not rely at all on the so-called packing lemma of Haussler \cite{Ha95}. Instead,
our approach is inspired by an argument used by Sudakov and Tomon \cite{ST20} in a related, but different
context.
It is perhaps important to also mention that, like Theorem \ref{FPSSZ1}, an asymmetric version of Theorem \ref{main} is
also very likely to have multiple applications in incidence geometry (with bounds which are superior
to the ones implied by Theorem \ref{FPSSZ1}). We will however not pursue such applications in this paper.

\bigskip

\section{Proof of Theorem \ref{main}}

\bigskip

Our notation is mostly standard. For a graph $G$ and vertex $v\in V(G)$, we write $N_G(v)$ for the set of neighbors of $v$ in $G$. When the graph is clear, we often write $N(v)$ for the same set. If $T\subset V(G)$, we write $N(T)$ for the set of common neighbours of the set $T$.

Let $G=(A,B,E)$ be a $K_{k,k}$-free bipartite graph with the number of edges satisfying $|E(G)| \geq cn^{2-1/d}$ for some constant $c > 0$, and where $|A|=|B|=n/2$. In order to show that $G$ has VC-dimension at least $d+1$, we need to prove the existence of a set $S$ of $d+1$ vertices in $A$, which is {\it{shattered}} by the set system formed by the neighborhoods $N(b)$ of the vertices $b \in B$; that is, a set $S$ with $|S|=d+1$ such that for every subset $S' \subset S$, there exist $b \in B$ with the property that $N(b) \cap S = S'$. 

First, we move to a subgraph with large minimum degree and choose the vertex which will be a neighbor to every vertex in $S$.

\begin{proposition} \label{lemma:subgraph}
Let $G$ be a $K_{k,k}$-free bipartite graph with parts $A$ and $B$ and at least $cn^{2-1/d}$ edges, where $c > 0$ is a constant and $|A|,|B|= n/2$. Then $G$ has an induced subgraph $G'$ with parts $A'\subset A$ and $B'\subset B$ such that $G'$ has minimum degree at least $\frac{c}{4}n^{1-1/d}$ and there exists a vertex $x\in B'$ such that $|N(x')\cap N(x)|=o(|N(x)|)$ holds for every $x'\in B'\setminus \{x\}$. 
\end{proposition}



In the short proof we shall use the asymmetric K\H ov\'ari--S\'os--Tur\'an theorem.

\begin{lemma} \label{lemma:asymmetric}
Let $G$ be a $K_{k,k}$-free bipartite graph with parts of size $m$ and $n$. Then $G$ has at most $O_k(nm^{1-1/k}+m)$ edges.
\end{lemma}

\begin{proofof}{Proposition \ref{lemma:subgraph}}
First, by iteratively discarding vertices of degree less than $\frac{c}{2}n^{1-1/d}$, we find a non-empty subgraph $G''$ which has minimum degree at least $\frac{c}{2}n^{1-1/d}$. Let $G''$ have parts $A''$ and $B''$. Choose a vertex $x\in B''$ arbitrarily. By Lemma \ref{lemma:asymmetric}, since $G''$ is $K_{k,k}$-free and $|N_{G''}(x)|\geq \frac{c}{2}n^{1-1/d}$, it is easy to see that the number of vertices $y\in B''$ such that $|N_{G''}(y)\cap N_{G''}(x)|\geq \frac{|N_{G''}(x)|}{\log n}$ is at most $n^{1/100}$. Write $B'$ for the set obtained from $B''$ after removing these vertices (apart from $x$). Let $A'=A''$ and let $G'$ be the induced subgraph of $G''$ with parts $A'$ and $B'$. This choice of $x$ and $G'$ satisfies the conditions in the statement of the proposition.
\end{proofof}

The next proposition will be applied for the subgraph found in the previous result.
 
\begin{proposition} \label{lemma:twocases}
Let $G$ be a bipartite graph with parts $A$ and $B$ and with minimum degree satisfying $\delta(G)=\delta \geq cn^{1-1/d}$ for some constant $c > 0$, and where $|A|,|B|\leq n/2$. Let $r$ be a constant positive integer and let $x\in B$. 
Then one of the following two statements must be true:

\begin{enumerate}
    \item there exists a set $R\subset N(x)$ of size $r$ such that for every $T\subset R$ of size $d$, we have $|N(T)|\geq r$ or
    \item there exist $\Theta(|N(x)|^r)$ sets $R\subset N(x)$ of size $r$ such that for every $T\subset R$ of size $d$, we have $N(T)\setminus \{x\}\neq\emptyset$.
\end{enumerate}
\end{proposition}

To prove this, we will use the so-called {\it{hypergraph removal lemma}}, proved independently by Nagle, R\"odl, Schacht \cite{NRS06} and Gowers \cite{G07}.

\begin{lemma} \label{lemma:hrl}
Let $r,d$ be positive integers. For every $\beta > 0$ there exists $\delta = \delta(r,d,\beta) > 0$ such that the following holds. If $\mathcal{H}$ is a $d$-uniform hypergraph on $N$ vertices such that one needs to remove at least $\beta N^{d}$ hyperedges of $\mathcal{H}$ in order to make it free of copies of $K_{r}^{(d)}$, then $\mathcal{H}$ contains at least $\delta N^{r}$ copies of $K_{r}^{(d)}$.
\end{lemma}

\begin{proofof}{Proposition \ref{lemma:twocases}}
Define a $d$-uniform hypergraph $\mathcal{H}$ on vertex set $N(x)$ such that any $T\subset N(x)$ of size $d$ is a hyperedge in $\mathcal{H}$ if and only if $N(T)\setminus \{x\}\neq \emptyset$. Then outcome 2. is equivalent to saying that $\mathcal{H}$ contains $\Theta(|N(x)|^r)$ copies of $K_r^{(d)}$. If the first statement does not hold, by using Lemma \ref{lemma:hrl}, it then suffices to prove that in order to destroy all copies of $K_r^{(d)}$ in $\mathcal{H}$, one needs to remove $\Theta(|N(x)|^d)$ hyperedges from $\mathcal{H}$. To this end, we shall prove that this indeed holds provided that there is no set $R\subset N(x)$ of size $r$ such that for every $T\subset R$ of size $d$, we have $|N(T)|\geq r$.

Color a $d$-set $T\subset N(x)$ green if $N(T)=\{x\}$, blue if $1<|N(T)|<r$ and red if $|N(T)|\geq r$. Note that if $y\in B\setminus \{x\}$, then any $d$-set $T\subset N(x)\cap N(y)$ is colored blue or red. If there exists some $R\subset N(x)$ of size $r$ such that every $T\subset R$ of size $d$ is red, then condition 1. holds. However, if $\ell$ is sufficiently large, by Ramsey's theorem (see \cite{R30} or \cite[Theorem 4.18]{Ju11}) we know that for every set $L\subset N(x)\cap N(y)$ of size $\ell$, there exists a subset $R\subset L$ of size $r$ such that all $d$-sets in $R$ have the same color.

Therefore each $\ell$-set in $N(x)\cap N(y)$ contains a monochromatic blue $r$-set. Clearly any such $r$-set $R$ has the property that for every $T\subset R$ of size $d$, we have $N(T)\setminus \{x\}\neq \emptyset$. Hence, if we are to delete all such $r$-sets from $\mathcal{H}$, then we need to delete a blue edge from every $\ell$-set in $N(x)\cap N(y)$, for every $y\in B\setminus \{x\}$. Hence, we need to delete at least $$\frac{1}{r}\sum_{y\in B\setminus \{x\}} \frac{1}{\binom{\ell}{d}} \binom{d(x,y)}{d}$$ hyperedges, where $d(x,y)=|N(x)\cap N(y)|$. Clearly, 
$$\sum_{y\in B\setminus \{x\}} d(x,y)=e(B\setminus \{x\},N(x))\geq \sum_{z\in N(x)} (d(z)-1)\geq d(x)(\delta-1).$$ Hence, by Jensen's inequality,
$$\sum_{y\in B\setminus \{x\}} \binom{d(x,y)}{d}\geq \Omega(n(d(x) \delta/n)^d)\geq \Omega(d(x)^d)=\Omega(|N(x)|^d).$$
Note that in the first inequality we have used implicitly that $d(x)\delta\geq (cn^{1-1/d})^2\geq \omega(n)$ as $d\geq 3$.

Thus, we indeed need to delete $\Omega(|N(x)|^d)$ hyperedges to destroy all copies of $K_r^{(d)}$ in $\mathcal{H}$.
\end{proofof}

\begin{proposition} \label{lemmma:ifallbig}
Let $q$ be a positive integer and let $F$ be a $K_{k,k}$-free bipartite graph with parts $B$ and $Q$, where $|Q|=q$. Assume that there exists $x\in B$ which is joined to all vertices of $Q$ and that for every $T\subset Q$ of size $d$, we have $|N(T)|\geq q$. If $q$ is sufficiently large compared to $d$ and $k$, then $Q$ has a subset of size $d+1$ that is shattered by $\{N(b):b\in B\}$.
\end{proposition}

\begin{proofof} {Proposition \ref{lemmma:ifallbig}}
Let $Z$ be a uniformly random subset of $Q$ of size $d+1$. We shall prove that $Z$ is shattered with probability at least $1/2$. 

In order to do this, we show that with probability at least $1/2$, the following property holds.

\begin{center}
{\it{For every $S\subset Z$ of size at most $d$ and every $z\in Z\setminus S$, we have $|N(S\cup \{z\})|< \frac{1}{d+1}|N(S)|$.}}
\end{center}

For convenience, let us call this property by the name {\it{Property $\mathcal{VC}$}} and see first why it implies that the set $Z$ is shattered. For each $S\subset Z$, we need to choose a vertex $b_S\in B$ such that $N(b_S)\cap Z=S$. For $S=Z$, choose $b_S=x$. Let $S\subset Z$ be a set of size at most $d$. By Property $\mathcal{VC}$, the number of vertices in $N(S)$ which have a neighbor in the set $Z\setminus S$ is less than $(d+1)\frac{1}{d+1}|N(S)|=|N(S)|$. Hence, we can pick some $b_S\in N(S)$ with $N(b_S)\cap Z=S$.

It remains to prove that Property $\mathcal{VC}$ holds with probability at least $1/2$. Using the union bound and conditioning on the $d$-subsets of $Z$, it suffices to prove that for any $S\subset Q$ of size at most $d$, the probability that there exists $z\in Z\setminus S$ with $|N(S\cup \{z\})|\geq \frac{1}{d+1}|N(S)|$ is at most $\frac{1}{2\cdot 2^{d+1}}$. However, note that since $|N(S)|\geq q$, where $q$ is sufficiently large compared to $d$ and $k$, and the induced subgraph of $F$ with parts $Q$ and $N(S)$ is $K_{k,k}$-free, it follows by Lemma \ref{lemma:asymmetric} that the number of vertices $y\in Q\setminus S$ with $|N(S\cup \{y\})|\geq \frac{1}{d+1}|N(S)|$ is at most $f(d,k)$ for some function $f$. Hence, if $q$ is sufficiently large, then with probability more than $1-\frac{1}{2\cdot 2^{d+1}}$, the random subset $Z\subset Q$ avoids all these vertices.
\end{proofof}

\bigskip

We are now in a position to prove Theorem \ref{main}.

\begin{proofof}{Theorem \ref{main}}
Suppose, for contradiction, that $G$ has at least $cn^{2-1/d}$ edges for some constant $c>0$. Choose a subgraph $G'$ and a vertex $x$ as in Proposition \ref{lemma:subgraph}. In what follows, all neighborhoods are defined in $G'$. By Proposition \ref{lemma:twocases}, one of the following two must hold.

\begin{enumerate}
    \item there exists a set $R\subset N(x)$ of size $r$ such that for every $T\subset R$ of size $d$, we have $|N(T)|\geq r$, or
    \item there exist $\Theta(|N(x)|^r)$ sets $R\subset N(x)$ of size $r$ such that for every $T\subset R$ of size $d$, we have $N(T)\setminus \{x\}\neq\emptyset$.
\end{enumerate}

If condition 1. holds, then by Proposition \ref{lemmma:ifallbig}, $G'$ has VC-dimension at least $d+1$. This implies that $G$ also has VC-dimension at least $d+1$, which is a contradiction.

So we may assume that condition 2. holds. Let $R\subset N(x)$ be a set of size $r$ such that for every $T\subset R$ of size $d$, we have $N(T)\setminus \{x\}\neq\emptyset$.

Let $q$ be a constant which is sufficiently large compared to $d$ and $k$. Now if $r$ is sufficiently large, by Ramsey's theorem, there exists a set $Q\subset R$ of size $q$ such that either $|N(T)|\geq q$ for every $T\subset Q$ of size $d$, or $|N(T)|<q$ for every $T\subset Q$ of size $d$. In the former case, Proposition \ref{lemmma:ifallbig} shows that $G'$ has VC-dimension at least $d+1$, which is a contradiction. Hence, $|N(T)|<q$ for every $T\subset Q$ of size $d$.

Since we can start with $\Theta(|N(x)|^r)$ many possible sets $R$ to get a subset $Q\subset R$ as above, it follows that there exist $\Theta(|N(x)|^q)$ sets $Q\subset N(x)$ such that for every $T\subset Q$ of size $d$, we have $N(T)\setminus \{x\}\neq \emptyset$ and $|N(T)|<q$.

Let $Q$ be such a set. Assume that the sets $N(T)\setminus \{x\}$ are pairwise disjoint as $T$ ranges over all subsets of $Q$ of size $d$. Then we distinguish between two cases.

Case (a) is when for every $S\subset Q$ of size $d-1$ we have $|N(S)|\geq q$. In this case, if $q$ is sufficiently large compared to $d$ and $k$, then using an argument very similar to the one in Proposition \ref{lemmma:ifallbig}, $Q$ has a subset of size $d+1$ which is shattered. In fact, just as in the proof of that proposition, a random subset $Z$ of size $d+1$ is shattered with probability at least $1/2$. The only difference is that for every set $S\subset Z$ of size $d$, the vertex $b_S$ is chosen from the set $N(S)\setminus \{x\}$. Since these sets are disjoint, we get a different vertex for each $S$, and we have $N(b_S)\cap Z=S$. Sets of size at most $d-1$ are treated just like in Proposition \ref{lemmma:ifallbig}. Now $G'$ has VC-dimension at least $d+1$, which is a contradiction.

Otherwise, which we call case (b), there exists $S\subset Q$ of size $d-1$ such that $|N(S)|<q$. However, note that the number of $q$-sets for which case (b) can occur is $o(|N(x)|^q)$. Indeed, there are $O(|N(x)|^{d-1})$ ways to choose the set $S$ and there are less than $q$ ways to choose a common neighbor of the set $S$. But any element $z\in Q\setminus S$ has a neighbor in $N(S)\setminus \{x\}$ (since $N(S\cup \{z\})\setminus \{x\}\neq \emptyset$), which leaves $o(|N(x)|)$ choices for each of these vertices. Altogether, we get only $o(|N(x)|^q)$ possibilities for $Q$.

It follows that the number of $q$-sets $Q\subset N(x)$ such that for every $T\subset Q$ of size $d$ we have $|N(T)|<q$ and the sets $N(T)\setminus \{x\}$ are not pairwise disjoint is $\Theta(|N(x)|^q)$. We now show that this is impossible by upper bounding the number of such sets. Let us choose distinct sets $T,T'\subset Q$ of size $d$ such that $(N(T)\setminus \{x\})\cap (N(T')\setminus \{x\})\neq \emptyset$. Note that there are at most $|N(x)|^d$ ways to choose $T$, given such a choice there are at most $|N(T)|<q$ ways to choose an element from $(N(T)\setminus \{x\})\cap (N(T')\setminus \{x\})$, and given these there are at most $o(|N(x)|)$ ways to choose each vertex in $T'\setminus T$. Every vertex in $Q\setminus (T\cup T')$ can be chosen in at most $|N(x)|$ many ways, so altogether we only have $o(|N(x)|^q)$ possibilities for $Q$, which is a contradiction.
\end{proofof}

\section{Concluding remarks}

It is natural to wonder what is the best upper bound that one could hope for in Theorem \ref{main}. We would like to make the following conjecture, which is somewhat related to Conjecture 1.2 from \cite{CL19}.

\begin{conjecture} \label{conj}
Let $k$ and $d$ be integers such that $k \geq d \geq 3$. Let $G = (A,B,E)$ be a bipartite graph with $|A|=|B|=n/2$ with VC-dimension at most $d$. Then, if $G$ is $K_{k,k}$-free, there exists some $\epsilon = \epsilon(d) > 0$ such that
$$|E(G)| = O_{n \to \infty}\left(n^{2-1/d-\epsilon}\right).$$
\end{conjecture}

It seems that proving Conjecture \ref{conj} would require finding a way to circumvent our usage of the hypergraph removal lemma or finding a new way to use the hypergraph removal lemma in a hypergraph with bounded VC-dimension instead (for which the implicit quantitative bounds would be much better, see for example \cite[Theorem 1.3]{FPS19} which would imply such a result). 

At the other end of the spectrum, we have the following lower bound, which can be naturally obtained from the probabilistic method. 

\begin{proposition} \label{lower bound}
    There exists a positive absolute constant $C$ with the following property. Let $k$ and $d$ be sufficiently large integers with $k\geq 2d-5$. Then there exists $c>0$ such that for every $n$ there is a $K_{k,k}$-free bipartite graph $G$ on $n+n$ vertices with VC-dimension at most $d$ such that $|E(G)|\geq cn^{2-\frac{1}{d-2-C/d}}$.
\end{proposition}

The main idea is to use following well-known result, which comes from a standard deletion argument. See \cite[Theorem 2.26]{FS13} and the references therein for more information.

\begin{lemma} \label{lemma:deletion method}
    Let $\mathcal{H}$ be a finite family of graphs. Let $\gamma=\max\left\{\frac{v(H)-2}{e(H)-1}: H\in \mathcal{H}\right\}$. Then there exists a positive constant $c$ such that for every $n$ there is an $n$-vertex graph with at least $cn^{2-\gamma}$ edges which does not contain any $H\in \mathcal{H}$ as a subgraph.
\end{lemma}

\begin{proofof}{Proposition \ref{lower bound}}
Let $F$ be the bipartite graph with vertex sets $X$ and $Y$, where $|X|=d+1$, for every $S\subset X$ of size at least $d-2$ there is a unique $y\in Y$ with $N_F(y)=S$ and there are no other vertices in $Y$. Observe that if a bipartite graph $G$ has VC-dimension at least $d+1$, then it contains $F$ as a subgraph. Moreover, note that $$v(F)=(d+1)+\binom{d+1}{3}+\binom{d+1}{2}+\binom{d+1}{1}+\binom{d+1}{0}$$ and
$$e(F)=(d-2)\binom{d+1}{3}+(d-1)\binom{d+1}{2}+d\binom{d+1}{1}+(d+1)\binom{d+1}{0}.$$

It is not hard to see that there exists a positive absolute constant $C$ such that for sufficiently large $d$, we have 
$$\frac{v(H)-2}{e(H)-1}\leq \frac{1}{d-2-C/d}.$$ Also, 
$$\frac{v(K_{k,k})-2}{e(K_{k,k})-1}=\frac{2k-2}{k^2-1}=\frac{2}{k+1}\leq \frac{1}{d-2}\leq \frac{1}{d-2-C/d}.$$

Thus, by Lemma \ref{lemma:deletion method}, there is a positive constant $c$ (dependent on $d$ and $k$) such that for any $n$ there is a graph with $2n$ vertices and at least $2cn^{2-\frac{1}{d-2-C/d}}$ edges which contains neither $F$ nor $K_{k,k}$ as a subgraph. By randomly partitioning its vertex set into two sets of size $n$ and taking the corresponding bipartite subgraph, it is now a standard linearity of expectation computation that proves there exists a bipartite graph $G$ with $n+n$ vertices and at least $cn^{2-\frac{1}{d-2-C/d}}$ edges which contains neither $F$ nor $K_{k,k}$ as a subgraph. Since $F$ is not a subgraph of $G$, the VC-dimension of $G$ is at most $d$, so $G$ is suitable.
\end{proofof}

\end{document}